\documentclass[reqno,   a4paper, 11pt]{amsart}
\usepackage{graphics}
\usepackage{color}
\usepackage{a4wide}
\usepackage{amsmath}
\usepackage{amssymb}
\usepackage{amsthm}
\usepackage{amsxtra}
\usepackage{latexsym}

\newcommand{\Nat}{{\mathbb{N}}}
\newcommand{\Z}{{\mathbb{Z}}}

\newcommand{\R}{{\mathbb{R}}}
\newcommand{\Q}{{\mathbb{Q}}}
\def\setsuchas#1#2{\left\{\,{#1}\,\vrule\,{#2}\,\right\}}
\newcommand{\set}[1]{{\{#1\}}}
\newcommand{\divides}[2]{{#1 \left\lvert {#2} \right.}}

\newcommand{\PowerRing}[4]{{#1 \lbrack\lbrack {#2}_{#3},
   \dots, {#2}_{#4} \rbrack \rbrack} } 

\newcommand{\PPxn}[1]{\PowerRing{#1}{x}{1}{n}}

\newcommand{\KKxn}{\PPxn{K}}

\usepackage{enumerate}

\theoremstyle{definition}
\newtheorem{definition}{\color{yellow}Definition}[section]
\newtheorem{example}[definition]{\color{green}Example}

\theoremstyle{remark}

\theoremstyle{plain}
\newtheorem{lemma}[definition]{\color{blue}Lemma}
\newtheorem{proposition}[definition]{\color{blue}Proposition}
\newtheorem{prop}[definition]{\color{blue}Proposition}
\newtheorem{theorem}[definition]{\color{blue}Theorem}

\newtheorem{corollary}[definition]{Corollary}

\newcommand{\tfh}{\overline{\pi_H}}
\newcommand{\tfZh}{\overline{\pi_{Z(H)}}}

\newcommand{\ptopat}{topologically prime atomic}
\newcommand{\topfact}{topologically factorial}
\newcommand{\ptopfact}{topologically prime factorial}

\newcommand{\monunits}[1]{{#1}^{\times}}
\newcommand{\free}[1]{\mathcal{F}{\left({#1}\right)}}
\newcommand{\atoms}{\mathcal{A}}
\newcommand{\topprimes}{\mathcal{B}}

\newcommand{\red}[1]{{#1}_{\mathrm{red}}}
\newcommand{\trivred}[1]{{#1}}
\newcommand{\mgenerates}[1]{\left\lbrack{#1}\right\rbrack}
\newcommand{\mtgenerates}[1]{\left\lbrack{#1}\right\rbrack^*}
\newcommand{\pr}{\operatorname*{pr}}

\newcommand{\propE}{allows arbitrary decimation}
\newcommand{\propEp}{allows finite decimation}
\newcommand{\propS}{allows dissociation}

\newcommand{\assumption}[1]{\fbox{\color{red}{#1}}}
\begin{document}

\sloppy
\author{Jan Snellman}
\address{Jan Snellman\\
Department of Mathematics, Stockholm University\\
106 91 Stockholm, Sweden}

\email{jans@matematik.su.se}
\urladdr{http://www.matematik.su.se/{\textasciitilde}jans}
\date{\today}

\title{Factorisation in topological monoids}

\thanks{\color{magenta}{Previous versions of this manuscript were written while the
  author was a post-doc at Laboratoire GAGE, Ecole Polytechnique,
  Palaiseau, France, and at University of Wales, Bangor, UK.,
supported by grants from \textit{Svenska Institutet}
  and \textit{Kungliga Vetenskapsakademin}
  and by grant n. 231801F from \textit{Centre International des Etudiants et
  Stagiaires}.}}

\begin{abstract}
  The aim of this paper is sketch a theory of divisibility and
  factorisation in topological monoids, where finite products are
  replaced by convergent products. The algebraic case can then be
  viewed as the special case of \emph{discretely topologized}
  topological monoids.

  In particular, we define the \emph{topological factorisation
    monoid}, a generalisation of the factorisation monoid for
  algebraic monoids, and show that it is always \emph{topologically
    factorial}: any element can be uniquely written as a convergent
  product of irreducible elements. 
\end{abstract}
\maketitle
\sloppy

\begin{section}{A primer on factorisation in discrete monoids}
    In this section, we give some basic definitions on the
  divisibility and factorisation theory of algebraic monoids.
  For additional information  we refer to
  \cite{HK:FinGenMon,ChapGer:krull}

    An (algebraic) monoid \(H\) is a semi-group with a neutral element.
    In this paper,  \assumption{\(H\) is assumed
    to be abelian and 
    cancellative}. Unless otherwise stated, we write \(H\)
    multiplicatively and denote by \( 1 \in H\) its neutral
    element. However, the monoid \(\Nat = (\Nat,+)\) will be written
      additively, with \(0 \in \Nat\) the neutral element. We will also write
      the monoid \(\Nat^X\) (with \(X\) a set) additively, along with its
      sub-monoids.

    We denote the set of units in \(H\) by \(\monunits{H}\), and say
    that \(H\) is \emph{reduced} if \(\monunits{H} = \set{1}\). Since
    \(\monunits{H}\) is a subgroup of \(H\), we can form the factor
    monoid \(\red{H} = H/\monunits{H}\), which is reduced. More
    explicitly, \(\red{H}=H/\sim\) where \(a \sim b\) iff \(a = eb\)
    for  some \(e \in \monunits{H}\). By passing to \(\red{H}\) if
    necessary, we will in what follows assume that \assumption{\(H\)
      is reduced}. We denote the set of non-units by \(H^*\),
    i.e. \(H^* = H \setminus \set{1}\). Since \(H\) is assumed to be
    both cancellative and reduced, it follows that it is also
    \assumption{torsion-free}. 
    
    If \(a,b \in H\) then we say that  \(a\)
    divides \(b\), and write \(\divides{a}{b}\), if there exists a
    (necessarily unique) element \(c \in H\) such
    that \(ac=b\). 
    An element \(p \in H^*\) is said
    said to be irreducible if \(p=a_1 \cdots a_r\) with \(a_i \in H\) implies that
     \(p=a_j\) for some \(j\). The irreducible elements in
    \(H\) are called atoms; we denote by \(\atoms(H)\) the set of all
    atoms in \(H\). 

    We say that \(p \in H^*\) is
    prime if whenever \(p\) divides \(a_1 \dots a_r\), it divides some \(a_i\).
    Note that a prime element is always irreducible.

    A monoid \(H\) is \emph{atomic} if \(p \in H\) can be written as a
        finite product of atoms, and \emph{factorial} if this
        factorisation is unique. In a factorial monoid, irreducible elements are prime. Hence
      unique factorisation into atoms implies unique factorisation into
      primes. A factorisation into primes is always unique, so if
      every element in \(H\) can be written as a product of primes,
      then \(H\) is factorial.
\end{section}

\begin{section}{Infinite products in topological monoids}
  We assume standard notions of topology and topological semigroups,
  as used for instance in
  \cite{topsem, Nagata:Topology,Bourbaki:Topology}. 

  In what follows, \(H\) will be a a \emph{topological monoid}, that is, a
  monoid with a 
    topology on its underlying set such that the multiplication map
    \(H \times H \to H\) is continuous. Following \cite{topsem}, we
    assume that the topology on \assumption{\(H\) is Hausdorff}. Note that any
    algebraic  monoid is a topological  monoid when endowed with
    the discrete topology. We also assume that \assumption{\(H\) is abelian,
    cancellative and reduced.}

    To be able to talk about infinite products, we
    make the following definition (this is the same definition as is
    used in \cite{Bourbaki:Topology}, but differently phrased):
    
    \begin{definition}\label{def:convprod}
      A (possibly infinite)  product 
      \begin{equation}
        \label{eqn:infprod}
        g = \prod_{j \in J} f_j, \qquad f_j \in H, \quad J \text{ any
          set} 
      \end{equation}
      is \emph{convergent} iff the net \((\phi,\Delta)\) converges to
      \(g\), where \(\Delta\) is the directed set of all finite
      subsets of 
      \(J\), and for \(F \in \Delta\),  \(\phi(F)=\prod_{j \in F}
      f_j\).

      In detail: for every neighbourhood \(U\) of \(g\), there is a
      finite subset \(S \subset J\), such that for all finite subsets
      \(S \subset T \subset J\) we have that \(\prod_{j \in T} f_j \in
      U\).
    \end{definition}

  \begin{definition}\label{def:propE}
    We say that \(H\)  \emph{\propE}
    if whenever \(b = \prod_{s \in S} 
    e_s\) is convergent, and \(T \subset S\), then \(\prod_{s \in
      T} e_s\) is convergent. We say that \(H\) \propEp\ if
    \(\prod_{s \in  H} e_s\) is convergent whenever \(S\setminus
    T\) is finite.
  \end{definition}
  For instance, if \(H\) is complete, then \(H\)
  \propE; this follows as in \cite[III, \S 5.3, Proposition
  3]{Bourbaki:Topology} (which treats the case of complete groups).
  
  \begin{example}\label{ex:propE}
    Let \(\R^+\) denote the monoid of non-negative real numbers, with
    the usual topology, and addition as the operation (which we'll
    write additively). Let \(\Q^+\) denote the sub-monoid of
    non-negative rational numbers. Then \(\Q^+\) is reduced,
    cancellative, Hausdorff, but not complete, and it does not allow
    arbitrary decimation. To see this, consider the sum 
    \begin{equation}
      \label{eq:sumdec}
      \sum_{k=1}^\infty 2^{-k} = 1,
    \end{equation}
    which (in  \(\R^+\)) can be decimated to yield any real number \(x
    \in [0,1]\).
  \end{example}

  \assumption{Henceforth, we assume that \(H\) \propE}.  

    We can regard any infinite product of the form 
    \begin{equation}
      \label{eqn:eqn:infmultiprod}
      U = \prod_{i \in I} g_i^{\alpha(i)}, \qquad \alpha(i) \in \Nat
    \end{equation}
    as an instance of \eqref{eqn:infprod} by taking \(J=I\) and
    \(f_j = g_i^{\alpha_i}\), or by replacing \(g^{\alpha(i)}\) with
    \(\alpha(i)\) copies of \(g\) and enlarging the index set
    accordingly.

  \begin{lemma}\label{lemma:Egivesdivisors}
    Suppose that 
    \(x=\prod_{i \in I} e_i^{\alpha(i)}\) is a convergent product in \(H\). 
    Then if \(\alpha(i) \ge \beta(i)\) for 
    all \(i \in I\) then the element 
    \(y=\prod_{i \in I} e_i^{\beta(i)}\) is a divisor of \(x\).
  \end{lemma}

    \begin{lemma}\label{lemma:propE}
      If \(x \in H^*\) then
      the sequence \((x^n)_{n=0}^\infty\) diverges.
    \end{lemma}
    \begin{proof}
      Suppose that \(x^n \to a \in H\). Then \(x \cdot x^n \to x a\),
      but 
      we also have that \(x \cdot x^n \to a\), hence \(xa = a\). This is
      impossible since \(H\) is cancellative and \(x \neq 1\).
    \end{proof}

    \begin{lemma}
       Suppose that 
       \(x=\prod_{i \in I} e_i^{\alpha(i)}\) is a convergent product
       in \(H\). Then all \(\alpha(i) < \infty\), and there is no
       infinite subset \(J \subset I\) such 
       that \(k, \ell \in J \implies e_k = e_\ell\).
    \end{lemma}
    \begin{proof}
      If \(\prod_{i \in I} e_i^{\alpha(i)}\) is convergent, then so is
      all its sub-products. By Lemma~\ref{lemma:propE}, this means that
      no element can occur infinitely many times in the product.
    \end{proof}

    \begin{theorem}\label{thm:multiset}
      Any convergent product \(A=\prod_{j \in J} e_j\) can be
      expressed as 
      \begin{equation}
        \label{eq:multiset}
        A = \prod_{h \in H^*} h^{m(h)}, \quad m(h)=\# \setsuchas{j \in
          J}{e_j = h} < \infty.
      \end{equation}
    \end{theorem}

    The study of infinite products in \(H\) is thus reduced to the study
    of certain multisets on \(H\). In the next section, we shall
    exploit a variant of this, when we consider multisets on the set
    of irreducible elements.

    \begin{definition}
      For \(M \subset H\) we 
      denote by \(\mtgenerates{M}\) the sub-monoid generated by all
      convergent products of elements in \(M\).
    \end{definition}

    \begin{lemma}\label{lemma:closeGen}
      If  \(M \subset H\), then 
      \begin{equation}
        \label{eqn:cls}
        M \subset \mgenerates{M} \subset         
        \mtgenerates{M} \subset 
        \overline{\mgenerates{M}} 
        \subset H
      \end{equation}
      where \(\overline{\mgenerates{M}}\) is the
      topological closure of \(\mgenerates{M}\).
    \end{lemma}
    \begin{proof}
      Note that the subspace \(\overline{\mgenerates{M}}\) is a
      sub-monoid of \(H\), since the closure of a sub-monoid is a
      sub-monoid. 
      
      The only non-trivial inclusion is 
      \( \mtgenerates{M} \subset  \overline{\mgenerates{M}}\). 
      Let, as in \eqref{eqn:infprod}, \(U = \prod_{j \in J} f_j\)
      be a convergent product, with \(J\) any  set; let \(\Delta\) be
      the 
      directed set of all finite subsets of  \(J\), let  for \(F \in
      \Delta\),  \(\phi(F)=\prod_{j \in F} f_j\), and let
      \((\phi,\Delta)\) be the corresponding net. By definition of
      convergent products, the net  converges to \(U \in H\). Since 
      \(\phi(F)=\prod_{j \in F} f_j \in \mgenerates{M}\) for all 
      \(F \in \Delta\) it follows that 
      \(U \in \overline{\mgenerates{M}}\). 
    \end{proof}
    
    \begin{example}\label{example:closenotstar}
      Let 
      \begin{displaymath}
        X = \set{0} \cup \setsuchas{1/n}{n \in \Nat^+} \subset
        [0,1] \subset \R        
      \end{displaymath}
      be given the subspace topology, and let \(H\)
      be the following abelian topological monoid. As an algebraic
      monoid, \(H\) is the free abelian monoid on \(X\); we denote by
      \(e_0\) the basis vector corresponding to 0, and by \(e_i\) the
      basis vector corresponding to \(1/i\). There is a
      surjection 
      \begin{align*}
        \phi: H & \to \R \\
        e_0 & \mapsto 0 \\
        e_i & \mapsto 1/i \qquad \text{ for } i > 0
      \end{align*}
      We give \(H\) the smallest topology such that \(\phi\) is
      continuous, that is, a sequence \(f_v \to f\) in \(H\) iff
      \(\phi(f_v) \to \phi(f)\) in \(\R\). 
      We will write the commutative monoid 
      operation on \(H\) additively.

      Let \(M = \setsuchas{e_n}{n \in \Nat^+} \subset H\). Since \(1/n
      \to 0\) in \(\R\), \(e_n \to e_0\) in \(H\), hence 
      \(e_0 \in    \overline{\mgenerates{M}}\). We claim that 
      \(e_0 \not \in \mtgenerates{M}\). To see this, let \(g: X \to \R\)
      be the natural  inclusion, which is of course continuous and
      closed. Suppose  that \(e_0 = \sum_{i \in I} e_{c_i}\) with 
      \(c_i \in \Nat^+\) is a convergent sum, 
      then by continuity of \(g\) we get that 
      \(0 = g(e_0) = g(\sum_{i \in I} e_{c_i}) =\sum_{i \in I}
      g(e_{c_i})\). 
      But \(g(e_{c_i}) > 0\) in the
      natural total order on \(\R\), hence 
      \(\sum_{i \in I} g(e_{c_i}) >   0\), a contradiction.
    \end{example}

    \begin{definition}\label{def:almostdiscrete}
      We say that the topological monoid \(H\) is \emph{almost
        discrete} if all convergent products in \(H\) are finite, that
      is, for all \(M \subset H\), we  have that 
      \(\mgenerates{M} = \mtgenerates{M}\).
    \end{definition}

    \begin{example}
      Let \(M\) be the multiplicative monoid of
      \(\KKxn\), with the
      \(\mathbf{m}\)-adic topology, and let \(H = \set{1} \cup {f \in
        M}{f(0)=0}\).  Then \(H\) is almost discrete (any convergent
      infinite product of elements in \(\KKxn\) will converge to \(0 \not \in
      H\), but not discrete (every non-polynomial in \(H\) is the limit of
      polynomials).  
    \end{example}

\end{section}

\begin{section}{Factorisation: atoms, primes, and their topological counterparts}
  \begin{definition}
    We say that \(p \in H\) is \emph{topologically prime} if whenever
    \(p\) divides a convergent product, it divides some factor. 
    We denote by \(\topprimes(H)\)
    the set of topologically prime elements in \(H\).
    \end{definition}

    \begin{lemma}\label{lemma:topirr}
    Call \(a \in H^*\) \emph{topologically
      irreducible} if it can not 
    be written as a convergent product of elements, all 
    different from \(a\). Then \(a\) is topologically irreducible if
    and only if it is irreducible.
    \end{lemma}
    \begin{proof}
      Since \(H\) is reduced, \(a\) is irreducible if and only if it
      can not be written as a finite product of non-units, all 
    different from \(a\). Thus if \(a\) is topologically irreducible,
    it is irreducible. 

    For the converse, we need to use that \(H\) \propEp. Suppose that
    \(a\) is irreducible, and that \(a = \prod_{i \in I} e_i\). Let
    \(j \in I\). Then \(a = e_j b\), where \(b=\prod_{i \in I
      \setminus \set{j}} e_i\). Since \(a\) is irreducible, either
    \(a= e_j\) or \(a=b\). If \(a=e_j\), we are done. If \(a=b\), then
    \(a = e_j b = b\), which is impossible since \(H\) is cancellative.
    \end{proof}

    \begin{definition}
      Suppose that 
      \begin{equation}\label{eq:fa}
        f=\prod_{i \in I} a_i = \prod_{j \in J} b_j
      \end{equation}
      are two convergent factorisations of \(f\) into non-units. We
      say that these factorisations are equivalent if 
      \begin{equation}
        \label{eq:equality}
        \forall h \in H^*: \qquad \# \setsuchas{i \in I}{a_i = h} \quad
        = \quad \# \setsuchas{i \in I}{a_i = h}
      \end{equation}
    \end{definition}

    \begin{definition}
      \(H\) is \emph{topologically} [\emph{atomic}, \emph{prime atomic}] if every
      non-unit \(h\) can be written as a convergent product of [atoms,
      topologically prime elements]. It is \emph{topologically}
      [\emph{factorial, prime factorial}] if any two such
      factorisations of \(h\) are equivalent.
    \end{definition}

    \begin{proposition}\label{prop:topFreePrime}
      Suppose that \(H\) is \ptopat, and let \(x \in H^*\). 
      Then the following are
      equivalent: 
      \begin{enumerate}[(i)]
      \item \(x\) is topologically prime,
      \item \(x\) is prime,
      \item \(x\) is irreducible.
      \end{enumerate}
    \end{proposition}
    \begin{proof}
      It suffices to show that an irreducible element is topologically
      prime, so suppose  that \(x\) is irreducible. Since \(H\) is \ptopat,
      \(x\) may be written as a convergent product
      of topologically prime elements. We claim that this product must
      have only one  
      factor. Hence, \(x\) is topologically
      prime. 

      To establish the claim, we argue by contradiction, and write \(x =
      \prod_{i \in I} e_i\) with \(e_i\) topologically prime, and where
      \(I\) is not a singleton. Choose an \(j \in I\) and put \(y=e_j\)
      and \(z=\prod_{i \in I \setminus \set{j}}e_i\). Since \(H\) \propEp, the
      latter product 
      is convergent. Hence \(x=yz\), in
      contradiction to the fact that \(x\) is irreducible.
    \end{proof}

    \begin{proposition}\label{prop:topfactPrime}
      Suppose that  \(H\) is  \topfact.
      Then atoms in \(H\) are prime.
      \end{proposition}
    \begin{proof}
      Suppose that \(x\) is an atom in \(H\), and that \(\divides{x}{ab}\),
      with \(a,b \in H\). Then there exists \(y \in H\) with \(xy=ab\).
      Since \(H\) is \topfact, we can uniquely factor \(a,b,y\) into atoms:
      \begin{displaymath}
        a = \prod_{i \in I} u_i, \quad 
        b = \prod_{j \in J} v_j, \quad 
        y = \prod_{k \in K} w_k.
      \end{displaymath}
      We can assume that \(I,J,K\) are pair-wise disjoint. Applying
      \cite[III, \S 5.3, Proposition
      3]{Bourbaki:Topology}\footnote{That Proposition deals with
        topological groups, but the proof works verbatim for
        topological monoids.} we have that
      \begin{displaymath}
        x \prod_{k \in K} w_k = \prod_{i \in I} u_i \, \times 
         \prod_{k \in K} w_k  = \prod_{\ell \in I \cup J} z_\ell, 
        \quad z_\ell = 
        \begin{cases}
          u_i & \text{ if } \ell = i \in I, \\
          v_j & \text{ if } \ell = j \in J. 
        \end{cases}
      \end{displaymath}
      Since \(H\) is \topfact, factorisation into atoms is unique, hence \(x = 
      z_\ell\) for some \(\ell \in I \cup J\). Without loss of generality,
      assume that \(\ell = i \in I\), so that \(x = u_i\). Then the fact that
      \(H\) \propEp\ implies that \(\divides{x}{a}\).
    \end{proof}
    
    \begin{definition}\label{def:propS}
      \(H\) \emph{\propS} if
    whenever \(b = \prod_{s \in S}
    e_s\) is convergent, and for each \(s \in S\), \(e_s =
    \prod_{l \in G_s} f_l\) is convergent, then 
\(b = \prod_{l \in G}  f_l\) is convergent, where \(G\) is the
      disjoint union of the \(G_s\)'s.
    \end{definition}
    
    \begin{example}
      There are lots of non-reduced, but cancellative and even
      complete monoids which
      do not allow expansion. For instance, in the additive group of
      the reals, if we put \(e_k = k + (-k) = 0\) we have that 
      \begin{math}
        0=\sum_{k \in \Z^+} (k+(-k)), 
      \end{math}
      but \(\sum_{k \in \Z} k\) is not convergent.
    \end{example}

    \begin{example}\label{ref:propS}
      Consider the monoid \(\Q^+\) of Example~\ref{ex:propE}. We claim
      that this monoid does \propS, but does not allow us to perform
      the reordering 
      \(\sum_{i \in I} \sum_{j \in J} e_{ij} = 
      \sum_{j \in J} \sum_{i \in I}  e_{ij}\).
      To establish the first claim, we note that in \(\R^+\) all
      summable nets are 
      countable, thus we need only to consider sequences. Furthermore,
      since everything is positive, all convergent sums in \(\R^+\)
      are absolutely convergent. Thus, if \(\sum_{i=1}^\infty
      \sum_{j=1}^\infty a_{ij} = A\), then \(\sum_{(i,j)}
      a_{ij}=A\). If the first sum has all summands in \(\Q^+\) and
      converges to a rational value, then so does the second sum.

      However, the ``column sums'' \(\sum_{i=1}^\infty a_{ij}\)
      need no be rational. 
    \end{example}
    
    We believe that our assumptions on \(H\) (cancellative, reduced)
    are \emph{not} enough to guarantee that \(H\) \propS, thus we
    postulate this in the next proposition.

\begin{proposition}\label{prop:topfactTopPrime}
      Suppose that  \(H\) is  \topfact\ and \propS.
      Then atoms in \(H\) are topologically prime.
      \end{proposition}
    \begin{proof}
      Suppose that \(x\) is an atom in \(H\). By the previous
      proposition, \(x\) is prime.
      Suppose that \(\divides{x}{\prod_{i \in I} f_i}\),
      with \(f_i \in H\), and write each \(f_i\) as a convergent
      product \(f_i = \prod_{a \in \atoms(H)} a^{g_i(a)}\).
      \begin{displaymath}
        \prod_{i \in I} f_i = \prod_{i \in I} \prod_{a \in \atoms(H)}
        a^{g_i(a)}, 
      \end{displaymath}
      but on the other hand,
      \(\prod_{i \in I} f_i = x b\) for some \(b \in H\).
      Write \(b = \prod_{a \in \atoms(H)} a^{h(a)}\), then
      \begin{displaymath}
        \prod_{i \in I} \prod_{a \in \atoms(H)} a^{g_i(a)} =
        x \prod_{a \in \atoms(H)} a^{h(a)}.
      \end{displaymath}
      Since \(H\) \propS, we get that
      \begin{displaymath}
        x \prod_{a \in \atoms(H)} a^{h(a)} = \prod_{(i,a) \in I \times
          \atoms(H)} a^{g_i(a)},
      \end{displaymath}
      and since \(H\) is \topfact, factorisation into atoms is unique,
      hence \(x\) occurs in the right hand side.
    \end{proof}

    \begin{proposition}\label{prop:Cuniq}
      If \(b = \prod_{s \in S} e_s^{\alpha_s}\) is a convergent
      product in \(H\) of topologically prime elements, then
      \(\alpha_s\) is the maximal integer \(r \ge 0\) such that
      \(\divides{e_s^{r}}{b}\).
      Hence, the factorisation of an element into
      a convergent product of topologically prime elements, if it
      exists, is unique.
    \end{proposition}
    \begin{proof}
      For any topologically prime element \(p \in M\), we have that 
      \(\divides{p}{b}\) iff \(p=e_s\) for some \(s \in S\). To see
      this, first note that if \(p=e_s\) then 
      \(b = e_s (e_s^{\alpha_s -1} \prod_{t \in S \setminus
        \set{s}}e_t^{\alpha_t})\).  The right hand side is a product of
      \(e_s\) and a convergent since \(H\) \propEp.

      Conversely, if
      \(\divides{p}{b}\) then by definition of topologically primeness
      there is an \(s \in S\) such that
      \(\divides{p}{e_s^{\alpha_s}}\). Applying the property of
      topologically primeness again, we get that \(\divides{p}{e_s}\).
      Clearly,   different topologically prime elements do not
      divide each other, hence \(\divides{p}{e_s}\) iff \(p=e_s\).
      It follows that
      \(\alpha_s\) is the maximal integer \(r\) such that
      \(\divides{e_s^{r}}{b}\). 
    \end{proof}

    \begin{corollary}\label{corr:topprimat}
      A \ptopat\ monoid  is \ptopfact.
    \end{corollary}

    So, we have the following implications: 
    \begin{displaymath}
      \text{\ptopat} \iff \text{\ptopfact} \implies \text{\topfact}.
    \end{displaymath}

    The last implication can be reversed if \(H\) \propS.

\end{section}

\begin{section}{The topological factorisation homomorphism}
  Recall that the free abelian monoid \(\free{\atoms(H)}\) is called
    the \emph{factorisation monoid} of \(H\), and the canonical
    homomorphism 
    \begin{equation}
      \label{eqn:pifacthom}
      \pi_H: \free{\atoms(H)}  \to H
    \end{equation}
    is called the \emph{factorisation homomorphism}. If \(p \in  H\), 
    then the elements in \(\pi_H^{-1}(p)\) are called the
    \emph{factorisations} of \(p\). This homomorphism gives a lot of
    information about the factorisation properties of \(H\): we have
    that \(H\) is atomic iff \(\pi_H\) is surjective, and factorial
    iff \(\pi_H\) is bijective.

    We now make a construction which captures also the infinite
    factorisations. First, we introduce some notation
    for the topological monoid \(\Nat^X\), the set of all functions
    \(X \to \Nat\), where \(X\) is any set. This is a topological
    monoid with the operation of point-wise addition (we'll write the
    operation additively), and the topology
    of point-wise convergence. It is also a partially ordered set with
    point-wise comparison. 

    \begin{definition}\label{def:Nadd}
      For \(x \in X\), we define
      \begin{align*}
        \chi_x: X & \to \Nat \\
        \chi_x(y) & = 
        \begin{cases}
          1 & x = y \\
          0 & x \neq y
        \end{cases}        
      \end{align*}
      Thus any \(f: X \to \Nat\) may be written as a  convergent
      sum \(f=\sum_{x \in X} f(x)   \chi_x\). 
    \end{definition}

    \begin{definition}
      The partially defined map
    \begin{align}
      \tfh: \Nat^{\atoms(H)} & \to H \\
      \label{eqn:toppifacthom}
      \sum_{a \in \atoms(H)} v(a) \chi_a & \mapsto \prod_{a \in \atoms(H)}
      a^{v(a)} 
    \end{align}
    is defined whenever the right-hand side of
    \eqref{eqn:toppifacthom} is a convergent product. Denote by
    \(Z(H) \subset \Nat^{\atoms(H)}\)
    the domain of definition of \(\tfh\).
    We call \(Z(H)\) \emph{topological factorisation monoid} of \(H\); note
    that it     contains the factorisation monoid of \(H\), since the
    latter corresponds to the finitely supported maps \(\atoms(H) \to \Nat\).
    In what follows, we regard \(\tfh\) as a map  \(Z(H) \to H\) and call it
    the \emph{topological factorisation homomorphism}.

    If \(p \in H\), then \(\tfh^{-1}(p)\) is the set of (topological)
    factorisations of \(p\).

    Clearly, \(H\) is topologically atomic iff \(\tfh\) is surjective, and
      topologically factorial iff \(\tfh\) is bijective.

    For each \(a \in \atoms(H)\), the projection map \(\pr_a:
    Z(H) \to \Nat\) is defined by \((e_v)_{v \in \mathcal{A(H)}}
    \mapsto e_a\).
    We topologize \(Z(H)\) by giving it the initial topology with respect to 
    \(\tfh\) and all the projection maps \(\pr_a\), where
    \(\Nat\) is 
    discretely topologized. 
    \end{definition}
    
    Thus, \(Z(H)\) has the weakest topology such that \(\tfh\) and
    all the 
    projections \(\pr_a\) are continuous, and a net \(f_i\) converges to \(f\)
    in \(Z(H)\) iff 
    \(\tfh(f_i) \to \tfh(f)\) and
    \(\forall a \in \mathcal{A}: \, \pr_a(f_i) \to \pr_a(f)\).
    It is easy to see that this topology is Hausdorff.

    \begin{lemma}
    With respect to the component-wise partial order  on
    \(\Nat^{\atoms(H)}\), \(Z(H)\) is an order ideal, i.e. if \(c,d:
    \atoms(H) \to \Nat\),  \(c \in Z(H)\), and \(d(a) \le c(a)\) for
    all \(a \in \atoms(H)\), then \(d \in Z(H)\).
    \end{lemma}
    \begin{proof}
     This is a direct
    consequence of our assumption that \(H\) \propE.
    \end{proof}

    \begin{prop}\label{prop:zhhom}
      \(Z(H)\) is a topological monoid, and \(\tfh\) is an
      homomorphism of 
      topological monoids.
    \end{prop}
    \begin{proof}
      To show that \(Z(H)\) is a algebraic monoid, we must show that if \(f,g
      \in Z(H)\) then \(h=f+g \in \Nat^{\atoms(H)}\) is in fact in \(Z(H)\).
      So, we must
      show that \(\prod_{a \in \atoms(H)}a^{h(a)}\) is convergent. Let \(W\)
      be a  
      neighbourhood of \(\tfh(h) \in H\). Since multiplication in \(H\) is
      continuous, there is a  neighbourhood \(U\) of \(\tfh(f)\) and a
      neighbourhood \(V\) of \(\tfh(g)\) such that \(UV \subset W\). Since 
      \(\prod_{a \in \atoms(H)}a^{f(a)}\) and
      \(\prod_{a \in \atoms(H)}a^{g(a)}\) are convergent, there is a finite
      \(S \subset \atoms(H)\) such that for any finite subset \(S \subset T
      \subset \atoms(H)\) we have that 
      \(\prod_{a \in T}a^{f(a)} \in U\) and
      \(\prod_{a \in T}a^{g(a)} \in V\), hence
      \(\prod_{a \in T}a^{h(a)} \in UV\subset W\). Thus, 
      \(\prod_{a \in \atoms(H)}a^{h(a)}\) is convergent, showing that \(Z(H)\) 
      is an algebraic monoid. The product obviously converges to \(\tfh(h)\),
      so \(\tfh\) is a homomorphism of algebraic monoids. By definition,
      \(\tfh\) is continuous.

      It remains to see that addition in \(Z(H)\) is continuous. Let
      \(f_i \to f\), \(g_i \to g\) in \(Z(H)\), then by definition \(\tfh(f_i) 
      \to \tfh(f)\) and \(\pr_a(f_i) \to \pr_a(f)\), and likewise for
      \(g\). Since \(H\) is a topological monoid, \(\tfh(f_i + g_i) =
      \tfh(f_i)\tfh(g_i) \to \tfh(f)\tfh(g) = \tfh(f+g)\), and similarly
      \(\pr_a(f_i+g_i) \to \pr_a(f+g)\). Thus \(f_i + g_i \to f + g\).
    \end{proof}

    \begin{lemma}\label{lemma:ZHpropE}
      \(Z(H)\) is reduced, cancellative, and \propE.
    \end{lemma}
    \begin{proof}
      Since \(\Nat^X\) is reduced and cancellative for all \(X\), we
      have that \(Z(H)\) is a submonoid of a reduced, cancellative
      monoid, and hence it is reduced and cancellative.

      Let \(\sum_{i \in I} g_i = h\) be a convergent sum in \(Z(H)\),
      and let \(J \subset I\). We want to show that \(\sum_{i \in J}
      g_i\) is convergent, i.e. that the following two conditions
      hold:
      \begin{enumerate}
      \item \(\prod_{i \in J} \tfh(g_i)\) converges,
      \item for all \(a \in \atoms{H}\), \(\sum_{i \in J} g_i(a) < \infty\).
      \end{enumerate}
      The first property follows since we know that \(\prod_{i \in I}
      \tfh(g_i)\) converges, and that \(H\) \propE. The second
      property follows since we know that \(\sum_{i \in I} g_i(a) < \infty\).
    \end{proof}

    We now show that this definition generalises the discrete one.
    \begin{theorem}
      If \(H\) is discrete, then so is \(Z(H)\).
    \end{theorem}
    \begin{proof}
      The fact that \(H\) is discrete means that \(Z(H)\) consists
      precisely of the finitely supported maps \(\atoms(H) \to H\).
      Suppose that \(c_i \to c \in Z(H)\). We put \(x=\tfh(c)\), and
      note that since \(\tfh(c_i) \to \tfh(c) = x\), and since \(H\) is
      discrete, there is an \(N_1\) such that \(\tfh(c_i)=x\) for all
      \(i > N_1\). 

      We have that \(c\) is supported on a finite set \(A \subset
      \atoms(H)\), and for each \(a \in A\), there is an \(M_a\) such
      that \(\pr_a(c_i) = \pr_a(c)\) whenever \(i > M_a\). Thus there
      is an \(N_2\) such that \(\pr_a(c_i) = \pr_a(c)\) whenever \(i >
      N_2\) and \(a \in A\). 

      Suppose now that \(i > N_1, N_2\). Then 
      \begin{displaymath}
        \tfh(c_i) = \prod_{a \in A} a^{\pr_a(c_i)} \prod_{b \in
          \atoms(H) \setminus A}  b^{\pr_b(c_i)} 
        = \tfh(c) =  \prod_{a \in A} a^{\pr_a(c)} \prod_{b \in
          \atoms(H) \setminus A}  b^{\pr_b(c)} 
        = \prod_{a \in A} a^{\pr_a(c_i)}, 
      \end{displaymath}
      so 
      \begin{displaymath}
        \prod_{b \in \atoms(H) \setminus A}  b^{\pr_b(c_i)} = 1,
      \end{displaymath}
      which means that \(\pr_b(c_i)=0\) for all \(b \not \in A\).
      Hence \(c_i = c\) when \(i > N_1, N_2\), so all convergent
      nets are stationary after a finite number of steps, which means
      that \(Z(H)\) has the discrete topology.
    \end{proof}

    \begin{lemma}\label{lemma:atoms}
      The atoms (and the topologically prime elements) of \(\Z(H)\)
      are precisely the elements \(\setsuchas{\chi_a}{a \in \atoms(H)}\)
    \end{lemma}
    \begin{proof}
      If \(a \in \atoms(H)\) then \(\chi_a \in Z(H) \subset
      \Nat^{\atoms(H)}\).  Since \(\chi_a\) is irreducible in
      \(\Nat^{\atoms(H)}\), it is irreducible in  \(\Z(H)\). If \(f
      \in Z(H)\), then \(f = \sum_{a \in \atoms(H)} f(a) \chi_a\), so
      if \(f(a), f(b) \neq 0\) for \(a \neq b\), then
      \(f\) is not topologically irreducible, hence
      (Lemma~\ref{lemma:topirr}) not irreducible. We have thus shown
      that \(\atoms(Z(H)) = \setsuchas{\chi_a}{a \in \atoms(H)}\).
      
      By Lemma~\ref{lemma:ZHpropE} we have that all \(\chi_a\) are
      topologically prime. A topologically prime element is prime,
      hence irreducible, hence of the form \(\chi_a\).
    \end{proof}

    \begin{lemma}\label{lemma:topprimdivthe}
      If \(H\) is a \topfact, and if
      \(H \ni f =  \prod_{a \in \atoms(H)} a^{v(a)}\),  
      \(H \ni g =  \prod_{a \in \atoms(H)} a^{w(a)}\),  
      then \(\divides{f}{g}\) iff \(\forall a \in \atoms(H): v(a)
      \le w(a)\).
    \end{lemma}
    \begin{proof}
      The underlying algebraic monoid of \(H\) is isomorphic
      to \(Z(H) \subset \Nat^{\atoms(H)}\), and \(Z(H)\) is an order ideal.
    \end{proof}

    \begin{example}\label{example:topdisc}
      For a \topfact\ monoid the topological factorisation homomorphism is an
      isomorphism of algebraic monoids, but it need not be an isomorphism of
      topological monoids. As an example, if
      \(H\) is the topological monoid of
      Example~\ref{example:closenotstar} then
      the factorisation
      homomorphism \(\tfh\) maps \(a_0\) to \(e_0\) and \(a_n\)
      to \(e_n\). Now, consider the sequence \(f_i = a_i\) for \(i > 1\).
      We have that 
      \(\tfh(f_i) \to e_0 =  \tfh(f_0)\). 
      However, \(\pr_{a_0}(f_i) = 0\) for all   \(i > 0\), so  
      \(\pr_{a_0}(f_i) \to 0\), whereas 
      \(\pr_{a_0}(f_0) =  1\). 
      Thus, the inverse is not continuous.
    \end{example}

      Recall that
      for a discretely topologized monoid, the factorisation monoid
      \(\free{\atoms(\trivred{H})}\) is free abelian. Similarly:      
      \begin{theorem}\label{thm:ZHptopfact}
        \(Z(H)\) is \ptopfact. 
      \end{theorem}
      \begin{proof}
        Each element \(f \in Z(H)\) can be written uniquely
        as 
        \(f = \sum_{a \in \atoms(A)} {\chi_a}^{f(a)}\); this sum
        is convergent with respect to the topology of
        point-wise convergence, and by construction, its image under \(\tfh\)
        is also convergent.
      \end{proof}

      Moreover:
      \begin{theorem}\label{thm:ZZH}
        \(Z(Z(H)) \simeq Z(H)\) as topological monoids.
      \end{theorem}
      \begin{proof}
        We define 
        \begin{equation}
          \label{eq:Xi}
          \begin{split}
            \Xi: Z(H) & \to \Nat^{\atoms(Z(H))} \\
            \sum_{a \in \atoms(H)} f(a) a & \mapsto \sum_{a \in
              \atoms(H)} f(a) \chi_a
          \end{split}
        \end{equation}
        where we have used the fact that \(\atoms(Z(H)) =
        \setsuchas{\chi_a}{a \in \atoms(H)}\)
        (Lemma~\ref{lemma:atoms}). This map is obviously an injective
        homomorphism of algebraic monoids.
        Since
        
        \begin{displaymath}
          \sum_{a \in  \atoms(H)} f(a) \chi_a \in Z(Z(H))   
        \end{displaymath}
        if and only if
        \begin{displaymath}
          \sum_{a \in \atoms(H)} f(a) a  \text{ converges in } Z(H)             
        \end{displaymath}
        if and only if
        \begin{displaymath}
          \prod_{a \in \atoms(H)} a^{f(a)} \text{ converges in } H  
        \end{displaymath}
        if and only if
        \begin{displaymath}
          \sum_{a \in \atoms(H)} f(a) a \in Z(H)
        \end{displaymath}
        we get that the image of \(\Xi\) is exactly \(Z(Z(H))\),
        thus that \(Z(H) \simeq Z(Z(H))\) as algebraic monoids.
        Henceforth, we regard \(\Xi\) as a map \(\Xi: Z(H) \to
        Z(Z(H))\). A simple calculation shows that it is a section to
        the topological factorisation homomorphism \(\tfZh: Z(Z(H))
        \to Z(H)\), i.e. that \(\tfZh \circ \Xi: Z(H) \to Z(H)\) is
        the identity.

        It remains to show that \(\Xi\) is continuous, with continuous
        inverse. Let \((g_v)_{v \in V}\) be a net in \(Z(H)\), and let
        \(g \in Z(H)\), with 
        \begin{equation}
          \label{eq:gvg}
          g = \sum_{a \in \atoms(H)} r(a) a, \qquad 
          g_v = \sum_{a \in \atoms(H)} r_v(a) a
        \end{equation}
        Then
        \begin{displaymath}
          g_v \to g   \text{ in } Z(H) 
        \end{displaymath}
         if and only if
         \begin{displaymath}
           \tfh(g_v) \to \tfh(g) \text{ in } H  \text{ and } \forall b
           \in \atoms(H): \pr_b(g_v) \to \pr_b(g)  
         \end{displaymath}
         if and only if
         \begin{displaymath}
            \tfh\left(\sum_{a \in \atoms(H)} r_v(a) a\right)  \to
            \tfh\left(\sum_{a \in 
              \atoms(H)} r(a) a\right) \text{ in 
              } H   
         \end{displaymath}
         and
         \begin{displaymath}
            \forall b 
            \in \atoms(H): \pr_b\left(\sum_{a \in \atoms(H)} r_v(a)
              a\right)  \to
            \pr_b\left(\sum_{a \in \atoms(H)} r(a) a\right)  
         \end{displaymath}
         which holds if and only if
         \begin{displaymath}
            \prod_{a \in \atoms(H)} a^{r_v(a)} \to \prod_{a \in
              \atoms(H)} a^{r(a)} 
             \text{ and }  \forall b \in \atoms(H): r_v(b) \to r(b).
         \end{displaymath}

         On the other hand,
         \begin{displaymath}
          \Xi(g_v) \to \Xi(g) \text{ in } Z(Z(H))             
         \end{displaymath}
         if and only if
         \begin{displaymath}
          \tfZh(\Xi(g_v)) \to \tfZh(\Xi(g)) \text{ in } Z(H) \text{ and } \forall b
            \in \atoms(Z(H)): \pr_b(\Xi(g_v)) \to \pr_b(\Xi(g)),           
         \end{displaymath}
          if and only if 
         \begin{displaymath}
             g_v \to g \text{ in } Z(H) \text{ and } \forall b \in
             \atoms(H): r_v(b) \to r(b)
         \end{displaymath}

        Hence, 
        \begin{displaymath}
          g_v \to g \iff \Xi(g_v) \to \Xi(g),
        \end{displaymath}
        so  \(\Xi\) is continuous, with continuous
        inverse.
      \end{proof}
\end{section}

\begin{section}{The case of restricted decimation --- a counterexample}
  If \(H\) does not \propE, then few of the previous results holds.
  We'll explore this in an example.

Consider the submonoid \(H \subset \Nat^\Nat\)
consisting of those functions \(\Nat \to \Nat\) which are either
finitely supported, or \(\ge 1\) everywhere.
Then, the the functions \(\chi_j\), with \(j \in \Nat\),
are topologically prime, and topologically irreducible. Now consider the
function \(f\) which is constantly \(1\).
Note that \(f\) is irreducible, since it can not be decomposed as a
non-trivial finite product. In fact, \(f\) is prime: if \(f\) divides
\(\prod_{j=1}^n g_j\) then at least one \(g_j\) must be \(\ge 1\)
everywhere, and then \(f\) divides that \(g_j\). 
However, \(f\) can be written as a convergent
(infinite) product of atoms, so it is neither topologically
irreducible, nor topologically prime.

Note that
\(f = \prod_{i=0}^\infty \chi_i\), yet \(\chi_0\) is not a divisor of
\(f\), since 
\begin{displaymath}
  \frac{f}{\chi_0} = \prod_{i = 1}^\infty \chi_i \not \in H. 
\end{displaymath}
So \(Z(H)\) is not an order ideal. This is true even if we instead
define the topological factorisation monoid as a set of maps from the
topologically irreducible elements to \(\Nat\).

\end{section}
\bibliographystyle{hplain}
\bibliography{journals,articles,snellman}
\end{document}